\documentclass[11pt,reqno]{amsart}

\usepackage[T1]{fontenc}
\usepackage{lmodern}
\usepackage[utf8]{inputenc}
\usepackage{amsmath,amssymb,amsthm,mathtools}
\usepackage{mathrsfs}
\usepackage{enumitem}
\usepackage[expansion=false]{microtype}

\usepackage[letterpaper,
            textwidth=0.85\paperwidth,
            hcentering,
            top=1in,
            bottom=1in,
            headsep=18pt,
            footskip=32pt]{geometry}

\usepackage[colorlinks=true,linkcolor=blue,citecolor=blue,urlcolor=blue]{hyperref}

\usepackage{parskip}
\usepackage{ragged2e}
\RaggedRight
\usepackage{titlesec}
\titleformat{\section}
    {\normalfont\bfseries\Large}{\thesection}{0.75em}{}
\titleformat{\subsection}
    {\normalfont\bfseries\large}{\thesubsection}{0.75em}{}
\titleformat{\subsubsection}
    {\normalfont\bfseries}{\thesubsubsection}{0.75em}{}
\titlespacing*{\section}{0pt}{1.5\baselineskip}{0.75\baselineskip}
\titlespacing*{\subsection}{0pt}{1.25\baselineskip}{0.5\baselineskip}

\makeatletter
\def\l@section{\@tocline{1}{0pt}{1pc}{2.5pc}{}}
\makeatother

\usepackage{fancyhdr}
\newcommand{\papershorttitle}{On integers $n$ relatively prime to $\lfloor \alpha n^2\rfloor$}
\fancypagestyle{paperstyle}{%
    \fancyhf{}
    \fancyhead[C]{\normalfont\scshape\small\papershorttitle}
    \fancyfoot[C]{\normalfont\small\thepage}

}
\fancypagestyle{plain}{%
    \fancyhf{}
    \fancyhead[C]{\normalfont\scshape\small\papershorttitle}
    \fancyfoot[C]{\normalfont\small\thepage}

}
\fancypagestyle{firstpage}{%
    \fancyhf{}
    \fancyhead[C]{\normalfont\scshape\small\papershorttitle}
    \fancyfoot[C]{\normalfont\small\thepage}

}

\setlist[itemize]{leftmargin=2em}

\theoremstyle{plain}
\newtheorem{theorem}{Theorem}[section]

\newtheorem{lemma}[theorem]{Lemma}
\newtheorem{corollary}[theorem]{Corollary}

\theoremstyle{definition}

\numberwithin{equation}{section}

\title[On integers $n$ relatively prime to $\lfloor \alpha n^2\rfloor$]{On integers $n$ relatively prime to $\lfloor \alpha n^2\rfloor$}
\date{2 November 2026}
\begin{document}

%\vspace{0.5em}
%{\large Aahan Chatterjee\par}
%\vspace{0.25em}
%{\small Universit\'e de Montr\'eal\par
%\href{mailto:aahan.chatterjee@umontreal.ca}{\nolinkurl{aahan.chatterjee@umontreal.ca}}\par}

\bigskip
\begin{abstract}
    For $\alpha\in\mathbb{R}$ and $k\in\mathbb{N}$, let $\mathcal S_{\alpha,k}(x):=\#\{n\in\mathbb N:n\le x,\ \gcd(n,\lfloor\alpha n^k\rfloor)=1\}$. Roth conjectured, and Watson proved, that $\mathcal S_{\alpha,1}(x)\sim 6x/\pi^2$ for every irrational $\alpha$. Whether the same asymptotic holds for every fixed $k\ge2$ is a special case of a question of Bergelson and Richter, and was asked explicitly by Champagne. We answer this for $k=2$, with no Diophantine condition on $\alpha$: for every irrational $\alpha$, $\mathcal S_{\alpha,2}(x)\sim x/\zeta(2)$ as $x\to\infty$. More precisely, for every fixed $A>0$ we obtain an error term $O_A(x/(\log q(x))^A)$, where $q(x)$ is the largest continued-fraction convergent denominator of $\alpha$ not exceeding $x$. The proof sieves $\mathcal S_{\alpha,2}(x)$ over divisors, splitting at the parameter $R=q(x)^{1/8}$. Divisors larger than $R$ are handled by reducing the condition $d\mid\gcd(n,\lfloor\alpha n^2\rfloor)$ to M\"obius sums over arithmetic progressions and applying Zhan's uniform estimate for exponential sums twisted by $\mu(n)$ over short intervals; divisors at most $R$ are handled instead through the Erd\H{o}s--Tur\'an inequality together with Weyl-type bounds for quadratic exponential sums of the form $\sum e(\beta n^2)$.
\end{abstract}
\author{Aahan Chatterjee}
\address{D\'epartement de math\'ematiques et de statistique,
Universit\'e de Montr\'eal, Montr\'eal, Canada}
\email{aahan.chatterjee@umontreal.ca}
\maketitle
\medskip
%\noindent\textbf{2020 Mathematics Subject Classification.} 11N25; 11L07, 11N36, 11J70.

%\noindent\textbf{Keywords.} Coprimality, integer parts, M\"obius function in short intervals, Weyl sums, continued fractions.

\bigskip
\tableofcontents

\section{Introduction}\label{sec:intro}
Roth conjectured that, for a fixed irrational real number $\alpha$,
the number of positive integers $n\le x$ that are relatively prime
to $\lfloor\alpha n\rfloor$ is asymptotic to $6x/\pi^2$
(see \cite{Watson1953}). More generally, for $\alpha\in\mathbb R$,
$k\in\mathbb N$ and $x>0$, define
\begin{equation*}
    \mathcal S_{\alpha,k}(x)
    :=\#\{n\in\mathbb N:n\le x,\
    \gcd(n,\lfloor\alpha n^k\rfloor)=1\}.
\end{equation*}
Watson~\cite{Watson1953} proved Roth's conjecture:
for every irrational $\alpha$,
\begin{equation*}
    \mathcal S_{\alpha,1}(x)\sim\frac{6x}{\pi^2}.
\end{equation*}

Bergelson and Richter~\cite{BergelsonRichter2017} showed that
the set of $n$ with $\gcd(n,\lfloor f(n)\rfloor)=1$ has density
$6/\pi^2$ for a wide class of functions $f$ from a Hardy field.
They worked under the growth hypotheses:
\[\frac{f(t)}{\log{t}\log\log\log\log{t}}\rightarrow \infty\qquad\text{as}\qquad t\rightarrow \infty\]
and that, for some
$j\in\mathbb N$,
\[
    \frac{f(t)}{t^{j-1}}\to\infty
    \qquad\text{and}\qquad
    f(t)=o(t^j)
    \qquad\text{as }t\to\infty.
\]
They also asked~\cite[Question~4]{BergelsonRichter2017} whether
this second hypothesis could be replaced by the conditions
\[
\begin{aligned}
    f(t)&=o(t^j)
        &&\text{for some }j\in\mathbb N\quad\text{and}\\
    \frac{|f(t)-p(t)|}{\log t}&\to\infty
        &&\text{for every }p\in\mathbb Q[t]\qquad\text{as }t\to\infty.
\end{aligned}
\]
For irrational $\alpha>0$ and $k\ge2$, the function
$f(t)=\alpha t^k$ satisfies these replacement conditions
but not the original growth hypothesis.

Under suitable Diophantine conditions on the $\alpha_j$ and for every
sequence of integers $1=m_1<m_2<\cdots<m_r$, Banks and
Shparlinski~\cite{BanksShparlinski2024} obtained asymptotic formulas with
power-saving error terms for the number of positive integers $n\le x$ such
that \[\gcd(n,\lfloor\alpha_1n^{m_1}\rfloor,\ldots,\lfloor\alpha_rn^{m_r}\rfloor)=1.\]
Champagne~\cite{Champagne2026} subsequently determined the density of such
$n$ for arbitrary irrational $\alpha_j$.
He also asked~\cite[Question~2]{Champagne2026} whether, for every
irrational $\alpha$ and every integer $k\ge2$, the set
\[
    \{n\in\mathbb N:\gcd(n,\lfloor\alpha n^k\rfloor)=1\}
\]
has density $1/\zeta(2)$, or equivalently, whether
$\mathcal S_{\alpha,k}(x)\sim6x/\pi^2$.

We answer this question for $k=2$, with no Diophantine
condition on $\alpha$:
\begin{theorem}\label{thm:main}
    Let $\alpha\in\mathbb R\setminus\mathbb Q$, and for $x\ge1$ let $q(x)$ denote the largest denominator of a continued-fraction convergent to $\alpha$ that does not exceed $x$. Then, for every fixed $A>0$ and all $x\ge x_0(\alpha)$,
    \begin{equation}\label{main_quantitative}
        \mathcal S_{\alpha,2}(x)
        =\frac{x}{\zeta(2)}
        +O_A\!\left(\frac{x}{(\log q(x))^A}\right),
    \end{equation}
    where the implied constant depends only on $A$, and $x_0(\alpha)$ depends only on $\alpha$. In particular, since $q(x)\to\infty$,
    \begin{equation}\label{conj_2}
        \mathcal S_{\alpha,2}(x)\sim\frac{x}{\zeta(2)}\qquad(x\to\infty).
    \end{equation}
\end{theorem}

\medskip\noindent\textbf{Overview of the proof.}
For a positive integer $d$ and $\alpha\in\mathbb R$, write $\mathcal A_{d,\alpha}(x):=\{n\in\mathbb N:n\le x, d\mid \gcd(n,\lfloor\alpha n^2\rfloor)\}$, so that the principle of inclusion--exclusion gives $\mathcal S_{\alpha,2}(x)=\sum_{d\le x}\mu(d)|\mathcal A_{d,\alpha}(x)|$. We split this sum at a parameter $R$ into a large-divisor piece $B_{\alpha,R}(x)$ ($d>R$) and a small-divisor piece $A_{\alpha,R}(x)$ ($d\le R$), and bound the two pieces by different methods.

For $B_{\alpha,R}(x)$ (Section~\ref{sec:large}), we reduce the condition $d\mid\lfloor\alpha n^2\rfloor$ to M\"obius sums over at most three pieces of arithmetic progressions, via the decomposition of Section~\ref{sec:decomp}, and bound the resulting sums using Zhan's uniform short-interval estimate (Section~\ref{sec:inputs}). For $A_{\alpha,R}(x)$ (Section~\ref{sec:small}), we instead use the equidistribution of $\{\alpha n^2\}$, via the Erd\H{o}s--Tur\'an inequality together with a Weyl-type bound for $\sum e(\beta n^2)$ (both recorded in Section~\ref{sec:inputs}; the latter is proved in Appendix~\ref{app:weyl}). Section~\ref{sec:proof} chooses $R=q(x)^{1/8}$, where $q(x)$ is the continued-fraction denominator of Lemma~\ref{uma_cont_frac}, to obtain an explicit error term and complete the proof.

The quadratic case provides a quantitative estimate without a Diophantine
restriction on $\alpha$. On the small-divisor side, the quadratic Weyl bound of
Lemma~\ref{uma_wooley} yields an error controlled by $q(x)$, even when $q(x)$ grows very slowly. In particular, the precise term
$N/\sqrt{s}$ is important for our argument
and even replacing this term by $N^{1+\varepsilon}/\sqrt{s}$, for any fixed
$\varepsilon>0$, would break the proof. So one of the obstacles to extending our methods to obtain a quantitative error term for $k\ge 3$ is that the standard Weyl bounds for sums of higher degree lose a factor depending on $N$, such as $N^{\varepsilon}$, which our argument cannot absorb when $q(x)$ grows slowly.

\section{Notation and Conventions}\label{sec:notation}
Throughout, $\alpha$ denotes a fixed real number, irrational unless stated otherwise. We write $e(t):=e^{2\pi i t}$ for $t\in\mathbb R$, $\{t\}$ for the fractional part of $t$, and $\|t\|:=\min_{n\in\mathbb Z}|t-n|$ for the distance from $t$ to the nearest integer. The letter $\mu$ denotes the M\"obius function, and $\mathsf 1_{E}$ the indicator function of a condition or set $E$. A subscript on $\ll$, $\gg$, or $\asymp$ records the permitted dependence of the implied constant; an implied constant carrying no subscript is absolute. All counting variables and divisors are positive integers. We use $\gcd(n,m)=\gcd(n,|m|)$ and $\gcd(n,0)=n$. For an interval $I$ we write $|I|$ for its length; intervals may be open, closed or half-open, and may be empty or consist of a single point. Table~\ref{tab:notation} collects the recurring pieces of notation used in the rest of the paper; each is (re)defined at its first occurrence in the text.

\begin{table}[ht]
\centering
\small
\begin{tabular}{ll}
\hline
Symbol & Meaning\\
\hline
$e(t)$ & $e^{2\pi i t}$\\
$\{t\}$ & fractional part of $t$\\
$\|t\|$ & distance from $t$ to the nearest integer\\
$|I|$ & length of the interval $I$\\
$\mu(n)$ & the M\"obius function\\
$\mathcal S_{\alpha,k}(x)$ & $\#\{n\in\mathbb N:n\le x,\gcd(n,\lfloor \alpha n^k\rfloor)=1\}$, Section~\ref{sec:intro}\\
$\mathcal A_{d,\alpha}(x)$ & $\{n\in\mathbb N:n\le x, d\mid \gcd(n,\lfloor\alpha n^2\rfloor)\}$, Section~\ref{sec:large}\\
$A_{\alpha,R}(x),\,B_{\alpha,R}(x)$ & small- and large-divisor parts of the sieve for $\mathcal S_{\alpha,2}(x)$, Section~\ref{sec:large}\\
$a(t),\,q(t)$ & continued-fraction convergent numerator/denominator at level $t$, Lemma~\ref{uma_cont_frac}\\
\hline
\end{tabular}
\caption{Recurring notation.}
\label{tab:notation}
\end{table}

\section{Some Standard Inputs}\label{sec:inputs}
In this section we collect the classical estimates used repeatedly in the rest of the paper: Zhan's estimate for exponential sums twisted by the M\"obius function over short intervals (used in Section~\ref{sec:decomp}), the Erd\H{o}s--Tur\'an discrepancy inequality and a Weyl-type bound for quadratic exponential sums (used in Section~\ref{sec:small}), and a quantitative form of Dirichlet's approximation theorem via continued fractions (used throughout).

\begin{theorem}[{\cite[Theorem~1]{Zhan1991}}]\label{mu_orth}
    For every fixed $A,\varepsilon>0$, uniformly for $\theta\in\mathbb R$
    and $2\le y\le u$ with $y\ge u^{2/3+\varepsilon}$, we have
    \begin{equation}\label{uma_zhan}
        \left|\sum_{u<n<u+y}\mu(n)e(n\theta)\right|
        \ll_{A,\varepsilon}\frac{y}{(\log y)^A}.
    \end{equation}
\end{theorem}
Zhan also extended the range to $y\ge u^{5/8+\varepsilon}$~\cite{Zhan1991b}; the range above suffices for our purposes.
As a direct consequence, we now record the following corollary to permit arbitrary interval lengths and endpoint
conventions.
\begin{corollary}\label{cor_zhan}
    Fix $A>0$ and $2/3<\sigma<1$. For every $D\ge2$, every interval
    $I\subseteq[D,2D]$ and every $\beta\in\mathbb R$, we have
    \begin{equation}\label{I}
        \left|\sum_{n\in I}\mu(n)e(n\beta)\right|
        \ll_{A,\sigma}\frac{|I|}{(\log D)^A}+D^\sigma.
    \end{equation}
    The implied constant is independent of $D$, $I$ and $\beta$.
\end{corollary}
\begin{proof}
    If $|I|<D^\sigma$, or if $D$ is bounded in terms of $\sigma$, the estimate follows from the trivial bound $\#(I\cap\mathbb Z)\le|I|+1\le D+1$.
    Otherwise, for sufficiently large $D$, we apply
    Theorem~\ref{mu_orth} with $u$ as the lower endpoint of $I$, $y=|I|$
    and $\varepsilon=(\sigma-2/3)/2.$ To see that the conditions of the theorem hold, note that
    \[
        |I|\ge D^\sigma
        \ge(2D)^{2/3+\varepsilon}
        \ge u^{2/3+\varepsilon},
        \qquad |I|\le D\le u.
    \]
    The final estimate now follows because $\log|I|\ge\sigma\log D$ and since the possible contribution of the endpoints $u$ and $u+y$ to the sum is $O(1).$
\end{proof}

\begin{theorem}[Erd\H{o}s--Tur\'an inequality; see {\cite[Ch.~1]{GranvilleRudnick2007}} or {\cite[Ch.~1]{Montgomery1994}}]\label{thm_et}
    Let $N,H\ge1$ be integers. For every choice of real numbers
    $\xi_1,\ldots,\xi_N$ and $0\le c_1<c_2\le1$, we have
    \begin{equation}\label{et}
    \begin{split}
        &\left|\#\{1\le n\le N:c_1\le\{\xi_n\}<c_2\}
        -N(c_2-c_1)\right|\\
        &\qquad\ll\frac{N}{H+1}
        +\sum_{h=1}^{H}\frac1h
        \left|\sum_{n=1}^{N}e(h\xi_n)\right|,
    \end{split}
    \end{equation}
    where the implied constant is absolute. In particular, the estimate
    is uniform in the sequence and in the interval endpoints.
\end{theorem}

\begin{lemma}[Quadratic Weyl bound (see {\cite[Lemma~1, p.~16]{Soundararajan2007}} or {\cite[Ch.~2]{GranvilleRudnick2007}})]\label{uma_wooley}
    Let $\beta\in\mathbb R$, $r\in\mathbb Z$ and $s\in\mathbb N$ satisfy
    \[
        s\ge2,\qquad (r,s)=1,\qquad
        \left|\beta-\frac rs\right|\le\frac1{s^2}.
    \]
    Then, for every real $N\ge1$,
    \begin{equation}\label{wooley}
        \left|\sum_{1\le n\le N}e(n^2\beta)\right|
        \ll\frac{N}{\sqrt{s}}+\sqrt{(N+s)\log s},
    \end{equation}
    where the implied constant is absolute.
\end{lemma}
This is a standard consequence of Weyl differencing. The precise term
$N/\sqrt{s}$ is important in Section~\ref{sec:small}, where we do not
assume a lower bound $s\gg N^\delta$ for any fixed $\delta>0$.
Replacing this term by $N^{1+\varepsilon}/\sqrt{s}$, for a fixed
$\varepsilon>0$, would therefore not, in general, suffice for our argument.
For completeness, we give a proof in Appendix~\ref{app:weyl}, along the
lines of the argument in \cite[Ch.~2]{GranvilleRudnick2007}.

\begin{lemma}\label{uma_cont_frac}
    Let $\beta\in\mathbb R\setminus\mathbb Q$. For each real $t\ge1$,
    let $a(t)/q(t)$, with $a(t)\in\mathbb Z$ and $q(t)\in\mathbb N$,
    be the last simple continued-fraction convergent to $\beta$ whose
    denominator does not exceed $t$.
    There is a constant $C(\beta)\ge16$ such that, for every real
    $t\ge C(\beta)$,
    \begin{equation}\label{uma_irr_cont}
        \left|\beta-\frac{a(t)}{q(t)}\right|\le\frac1{tq(t)},
        \qquad (a(t),q(t))=1,\qquad16\le q(t)\le t.
    \end{equation}
    Moreover, $q(t)\to\infty$ as $t\to\infty$.
\end{lemma}
\begin{proof}
    Let $\{a_i/q_i\}_{i\ge 0}$ be the simple continued-fraction convergents to $\beta$ and $j$ be the index so that $q_j=q(t)$ (that is, $q_j\le t<q_{j+1}$). Using the standard facts
    \[
        (a_i,q_i)=1,\qquad
        1=q_0\le q_1<q_2<\cdots,\qquad
        \left|\beta-\frac{a_i}{q_i}\right|
        <\frac1{q_iq_{i+1}}\qquad \text{for all}\quad i\ge 0
    \]
    (see \cite[Ch.~X]{HardyWright})
    gives
    \begin{equation}\label{uma_145}
        \left|\beta-\frac{a_j}{q_j}\right|
        <\frac1{q_jq_{j+1}}
        <\frac1{q_jt}
    \end{equation}
    and that $q_j\rightarrow \infty.$ So $q(t)\to\infty$.
    Thus taking $C(\beta)$ to be any convergent denominator
    at least $16$ ensures $q(t)\ge16$ for
    $t\ge C(\beta)$.
\end{proof}
\section{A Decomposition of the Condition
\texorpdfstring{$\{\beta d\}<1/d$}{\{beta d\}<1/d}}
\label{sec:decomp}

The large-divisor estimate of Section~\ref{sec:large} requires bounding
short sums of $\mu(d)$ restricted to the condition $\{\beta d\}<1/d$.
The purpose of this section is to bound these restricted short sums on any interval
$I\subseteq[D,2D].$ To do this, we will first show that the integers satisfying this condition form the union
of at most three pieces of arithmetic progressions with a common modulus
(Lemma~\ref{T_decomp}). We will then express these congruence conditions with
additive characters and then apply Corollary~\ref{cor_zhan} to each piece
separately (Lemma~\ref{umaze}).

Throughout this section, sums over $d\in I$ are over integers. Intervals
may have either endpoint included or excluded, and may be empty or
consist of a single point.

\begin{lemma}\label{umaze}
    Let $A>0$ and $2/3<\sigma<1$. Then for every $D\ge2$, every interval
    $I\subseteq[D,2D]$ and every $\beta\in\mathbb R$,
    \begin{equation}\label{target_ineq}
        \left|
        \sum_{d\in I}\mu(d)\mathsf{1}_{\{\beta d\}<1/d}
        \right|
        \ll_{A,\sigma}
        \frac{|I|}{(\log D)^A}+D^\sigma,
    \end{equation}
    where the implied constant depends only on $A$ and $\sigma$.
\end{lemma}

We prove Lemma~\ref{umaze} after establishing Lemma~\ref{T_decomp}.
Fix $D\ge1$, an interval $I\subseteq[D,2D]$, and
$\beta\in\mathbb R$, and put
\[
    \mathcal S_\beta
    :=\{d\in I\cap\mathbb Z:\{\beta d\}<1/d\}.
\]
By Dirichlet's approximation theorem
(see, e.g., \cite[Ch.~XI]{HardyWright}), there are coprime integers $a,q$
with $1\le q\le D$ such that
\begin{equation}\label{dirichlet}
    \left|\beta-\frac aq\right|<\frac1{qD}.
\end{equation}
The following decomposition lemma describes $\mathcal S_\beta$ in terms
of $a$ and $q$. It holds for every $q\ge1$: the small moduli $q\le4$ are
handled by the same case analysis as $q\ge5$ and need no additional
intervals.

\begin{lemma}\label{T_decomp}
    Let $D\ge1$, let $I\subseteq[D,2D]$ be an interval, let
    $\beta\in\mathbb R$, and let $a,q$ be coprime integers with
    $1\le q\le D$ satisfying \eqref{dirichlet}. Then there exist residues
    $r_1,r_2\pmod q$ and pairwise disjoint, possibly empty, subintervals
    $I_1,I_2,I_3\subseteq I$ such that
    \begin{equation}\label{target_decomp}
    \begin{split}
        \mathcal S_\beta
        &=
        \{d\in I_1\cap\mathbb Z:d\equiv r_1\pmod q\}\\
        &\qquad\cup
        \{d\in(I_2\cup I_3)\cap\mathbb Z:d\equiv r_2\pmod q\}.
    \end{split}
    \end{equation}
    Thus $\mathcal S_\beta$ is a disjoint union of at most three sets
    of the form $J\cap(r+q\mathbb Z)$, with $J\subseteq I$ an interval.
\end{lemma}

\begin{proof}
    Write $\beta=a/q+\gamma$, so that $|\gamma|<1/(qD)$ by
    \eqref{dirichlet}. Every $d\in I$ satisfies
    $q\le D\le d\le2D$, and therefore
    \begin{equation}\label{basic}
        |\gamma d|<\frac2q
        \qquad\text{and}\qquad
        \frac1d\le\frac1q.
    \end{equation}
    For $d\in I\cap\mathbb Z$, let $s\in\{0,1,\dots,q-1\}$ be the residue
    of $ad$ modulo $q$. Then $\beta d\in s/q+\gamma d+\mathbb Z$, and hence
    \begin{equation}\label{frac_s}
        \{\beta d\}=\left\{\frac sq+\gamma d\right\}.
    \end{equation}
    We determine $\mathcal S_\beta$ by a case analysis on $s$ and on the
    size of $\gamma d$. In each case, we use \eqref{frac_s} to identify
    $\{\beta d\}$ explicitly, and the condition $\{\beta d\}<1/d$ becomes
    a pair of simultaneous inequalities for $d$. We will show that only
    the following four systems occur:
    \begin{align}
        \gamma d&\ge0,
        &\gamma d&<\frac1d;
        \tag{$\mathrm T_0$}\label{T0}\\
        \gamma d&\ge\frac1q,
        &\gamma d-\frac1q&<\frac1d;
        \tag{$\mathrm T_1$}\label{T1}\\
        -\frac1q\le\gamma d&<0,
        &\gamma d+\frac1q&<\frac1d;
        \tag{$\mathrm T_{-1}$}\label{Tm1}\\
        \gamma d&<-\frac1q,
        &\gamma d+\frac2q&<\frac1d.
        \tag{$\mathrm T_{-2}$}\label{Tm2}
    \end{align}
    In each system $(\mathrm T_j)$ the first condition fixes the range of
    $\gamma d$, and the second is the condition $\{\beta d\}<1/d$ in the
    case $\{\beta d\}=\gamma d-j/q$. The non-strict lower bounds at the
    fractional-part boundary are intentional: values with
    $\{\beta d\}=0$ satisfy the target condition and must be retained.

    For each of the four systems, we first determine its real solution
    set in $I$ and show that it is a union of at most two intervals
    (Step~1). We then identify the residue classes and fractional-part
    ranges giving rise to these systems (Step~2), and combine the
    resulting sets to obtain the stated decomposition (Step~3).

    \smallskip\noindent
    \emph{Step 1: solving the four systems.}
    For each system we determine its solution set, that is, the set of
    real $t\in I$ satisfying it with $d$ replaced by $t$. Since
    $t\in I\subseteq[D,2D]$ and $D\ge1$, multiplying the second inequality
    by $t$ preserves its direction and turns it into a quadratic
    inequality.

    \eqref{T0}: If $\gamma<0$, the first inequality fails for every
    $t\in I$. If $\gamma=0$, both inequalities hold for every $t\in I$.
    If $\gamma>0$, the first inequality always holds, and the second is
    equivalent to $\gamma t^2<1$, that is, to $t<\gamma^{-1/2}$.
    So the solution set is $\emptyset$, $I$ or
    $I\cap(0,\gamma^{-1/2})$, an interval in each case.

    \eqref{T1}: If $\gamma\le0$, the first inequality has no solution
    in $I$. If $\gamma>0$, it is equivalent to
    $t\ge1/(q\gamma)$. The second becomes
    $\gamma t^2-t/q-1<0$. This quadratic has positive leading coefficient
    and its roots have product $-1/\gamma<0$, so for $t>0$ the inequality
    holds exactly when $t<w$, where
    \[
        w:=\frac{1/q+\sqrt{1/q^2+4\gamma}}{2\gamma}
    \]
    is the positive root. The solution set is
    $I\cap[1/(q\gamma),w)$, an interval.

    \eqref{Tm1}: If $\gamma\ge0$, the first inequality has no solution
    in $I$. If $\gamma<0$, it is equivalent to
    $t\le1/(q|\gamma|)$. The second becomes
    $\gamma t^2+t/q-1<0$, a quadratic inequality with negative leading
    coefficient. Put $\Delta:=1/q^2+4\gamma$. If $\Delta<0$, the quadratic
    is negative for every $t$, and the solution set is
    $I\cap(0,1/(q|\gamma|)]$. If $\Delta\ge0$, its roots are
    \[
        y_\pm:=\frac{1/q\pm\sqrt\Delta}{2|\gamma|},
        \qquad\text{which satisfy}\qquad
        0<y_-\le y_+\le\frac1{q|\gamma|},
    \]
    and it is negative exactly when $t<y_-$ or $t>y_+$. The solution set
    is then
    \[
        \bigl(I\cap(0,y_-)\bigr)
        \cup
        \bigl(I\cap(y_+,1/(q|\gamma|)]\bigr).
    \]
    When $\Delta=0$, the two roots coincide, and the common root is
    excluded because the quadratic inequality is strict. In either case
    the solution set is the union of at most two disjoint intervals.

    \eqref{Tm2}: If $\gamma\ge0$, the first inequality has no solution
    in $I$. If $\gamma<0$, it is equivalent to
    $t>1/(q|\gamma|)$. The second becomes
    $\gamma t^2+2t/q-1<0$. Since
    $|\gamma|<1/(qD)\le1/q^2$, we have
    $\Delta':=4/q^2+4\gamma>0$, so this quadratic has the two distinct
    roots
    \[
        z_\pm:=\frac{2/q\pm\sqrt{\Delta'}}{2|\gamma|},
    \]
    and it is negative exactly when $t<z_-$ or $t>z_+$. As
    $1/(q|\gamma|)$ is the midpoint of $z_-$ and $z_+$, we have
    $z_-<1/(q|\gamma|)<z_+$. Hence the solution set is
    $I\cap(z_+,\infty)$, an interval.

    \smallskip\noindent
    \emph{Step 2: the case analysis.}
    We treat $q\ge5$ first and then $q=4,3,2,1$. In each case we list the
    possible values of $s$, split according to the size of $\gamma d$,
    which by \eqref{basic} lies in $(-2/q,2/q)$, and identify
    $\{\beta d\}$ by \eqref{frac_s}. We say that a case
    \emph{gives nothing} if $\{\beta d\}\ge1/q$ throughout it; by
    \eqref{basic}, such a case contains no element of $\mathcal S_\beta$.
    Within a fixed residue class $ad\equiv s\pmod q$, we say that a case
    \emph{leads to} $(\mathrm T_j)$ if, under \eqref{basic}, its
    restrictions on $\gamma d$ together with $\{\beta d\}<1/d$
    are equivalent to that system.

    \emph{Case $q\ge5$.}
    If $d\in\mathcal S_\beta$, then by \eqref{basic}
    \[
    \begin{aligned}
        \left\|\frac{ad}{q}\right\|
        &=\|\beta d-\gamma d\|\\
        &\le\|\beta d\|+|\gamma d|\\
        &\le\{\beta d\}+|\gamma d|\\
        &<\frac1d+\frac2q\le\frac3q.
    \end{aligned}
    \]
    Since $\|ad/q\|=\min(s,q-s)/q$, it follows that
    $\min(s,q-s)<3$. As $s$ is an integer between $0$ and $q-1$, we obtain
    $s\in\{0,1,2,q-2,q-1\}$. These five residues are distinct since
    $q\ge5$.
    \begin{itemize}
        \item $s=0$: here $\{\beta d\}=\{\gamma d\}$. If $\gamma d<0$,
        then $\{\beta d\}=1+\gamma d>1-2/q\ge1/q$, so this case gives
        nothing. If $\gamma d\ge0$, then $\{\beta d\}=\gamma d$, and this
        case leads to \eqref{T0}.

        \item $s=1$: here $1/q+\gamma d\in(-1/q,3/q)$. If
        $\gamma d<-1/q$, then
        $\{\beta d\}=1+1/q+\gamma d>1-1/q\ge1/q$, so this case gives
        nothing. If $\gamma d\ge-1/q$, then
        $\{\beta d\}=\gamma d+1/q$; since $1/d\le1/q$, the inequality
        $\gamma d+1/q<1/d$ forces $\gamma d<0$, and this case leads to
        \eqref{Tm1}.

        \item $s=2$: here
        $2/q+\gamma d\in(0,4/q)\subseteq(0,1)$, so
        $\{\beta d\}=\gamma d+2/q$. Since $1/d\le1/q$, the inequality
        $\gamma d+2/q<1/d$ forces $\gamma d<-1/q$, and this case leads
        to \eqref{Tm2}.

        \item $s=q-1$: here
        $(q-1)/q+\gamma d\in((q-3)/q,(q+1)/q)$. If $\gamma d<1/q$,
        then $\{\beta d\}=(q-1)/q+\gamma d>1-3/q\ge1/q$, so this case
        gives nothing. If $\gamma d\ge1/q$, then
        $\{\beta d\}=\gamma d-1/q$, and this case leads to \eqref{T1}.

        \item $s=q-2$: here
        $(q-2)/q+\gamma d\in((q-4)/q,1)$, so
        $\{\beta d\}=(q-2)/q+\gamma d>1-4/q\ge1/q$, and this case
        gives nothing.
    \end{itemize}

    \emph{Case $q=4$.}
    Here $\gamma d\in(-1/2,1/2)$ and $s\in\{0,1,2,3\}$.
    \begin{itemize}
        \item $s=0$: if $\gamma d<0$, then
        $\{\beta d\}=1+\gamma d>1/2$, so this case gives nothing; if
        $\gamma d\ge0$, then $\{\beta d\}=\gamma d$, and this case leads
        to \eqref{T0}.

        \item $s=1$: here $1/4+\gamma d\in(-1/4,3/4)$. If
        $\gamma d<-1/4$, then $\{\beta d\}=5/4+\gamma d>3/4$, so this
        case gives nothing. If $\gamma d\ge-1/4$, then
        $\{\beta d\}=\gamma d+1/4$; the inequality
        $\gamma d+1/4<1/d\le1/4$ forces $\gamma d<0$, and this case
        leads to \eqref{Tm1}.

        \item $s=2$: here $1/2+\gamma d\in(0,1)$, so
        $\{\beta d\}=\gamma d+1/2$; the inequality
        $\gamma d+1/2<1/d\le1/4$ forces $\gamma d<-1/4$, and this case
        leads to \eqref{Tm2}.

        \item $s=3$: here $3/4+\gamma d\in(1/4,5/4)$. If
        $\gamma d<1/4$, then $\{\beta d\}=3/4+\gamma d>1/4$, so this
        case gives nothing. If $\gamma d\ge1/4$, then
        $\{\beta d\}=\gamma d-1/4$, and this case leads to \eqref{T1}.
    \end{itemize}

    \emph{Case $q=3$.}
    Here $\gamma d\in(-2/3,2/3)$ and $s\in\{0,1,2\}$.
    \begin{itemize}
        \item $s=0$: if $\gamma d<0$, then
        $\{\beta d\}=1+\gamma d>1/3$, so this case gives nothing; if
        $\gamma d\ge0$, then $\{\beta d\}=\gamma d$, and this case leads
        to \eqref{T0}.

        \item $s=1$: here $1/3+\gamma d\in(-1/3,1)$. If
        $\gamma d<-1/3$, then $\{\beta d\}=4/3+\gamma d>2/3$, so this
        case gives nothing. If $\gamma d\ge-1/3$, then
        $\{\beta d\}=\gamma d+1/3$; the inequality
        $\gamma d+1/3<1/d\le1/3$ forces $\gamma d<0$, and this case
        leads to \eqref{Tm1}.

        \item $s=2$: here $2/3+\gamma d\in(0,4/3)$. If
        $\gamma d<1/3$, then $\{\beta d\}=\gamma d+2/3$; the inequality
        $\gamma d+2/3<1/d\le1/3$ forces $\gamma d<-1/3$, and this case
        leads to \eqref{Tm2}. If $\gamma d\ge1/3$, then
        $\{\beta d\}=\gamma d-1/3$, and this case leads to \eqref{T1}.
    \end{itemize}

    \emph{Case $q=2$.}
    Here $\gamma d\in(-1,1)$ and $s\in\{0,1\}$.
    \begin{itemize}
        \item $s=0$: here $\{\beta d\}=\{\gamma d\}$. If $\gamma d<0$,
        then $\{\beta d\}=\gamma d+1$; the inequality
        $\gamma d+1<1/d\le1/2$ forces $\gamma d<-1/2$, and since
        $1=2/q$, this case leads to \eqref{Tm2}. If $\gamma d\ge0$,
        then $\{\beta d\}=\gamma d$, and this case leads to \eqref{T0}.

        \item $s=1$: here $1/2+\gamma d\in(-1/2,3/2)$. If
        $\gamma d<-1/2$, then $\{\beta d\}=3/2+\gamma d>1/2$, so this
        case gives nothing. If $-1/2\le\gamma d<1/2$, then
        $\{\beta d\}=\gamma d+1/2$; the inequality
        $\gamma d+1/2<1/d\le1/2$ forces $\gamma d<0$, and this case
        leads to \eqref{Tm1}. If $\gamma d\ge1/2$, then
        $\{\beta d\}=\gamma d-1/2$, and this case leads to \eqref{T1}.
    \end{itemize}

    \emph{Case $q=1$.}
    Here $\gamma d\in(-2,2)$, $s=0$, and
    $\{\beta d\}=\{\gamma d\}$.
    \begin{itemize}
        \item If $\gamma d<-1$, then $\{\beta d\}=\gamma d+2$, and this
        case leads to \eqref{Tm2}.

        \item If $-1\le\gamma d<0$, then $\{\beta d\}=\gamma d+1$, and
        this case leads to \eqref{Tm1}.

        \item If $0\le\gamma d<1$, then $\{\beta d\}=\gamma d$; as
        $\gamma d<1/d\le1$ implies $\gamma d<1$, this case leads to
        \eqref{T0}.

        \item If $\gamma d\ge1$, then $\{\beta d\}=\gamma d-1$, and this
        case leads to \eqref{T1}.
    \end{itemize}

    \smallskip\noindent
    \emph{Step 3: conclusion.}
    Step~2 shows that, for every integer $d\in I$, we have
    $d\in\mathcal S_\beta$ if and only if there exists
    $j\in\{0,1,-1,-2\}$ such that $d$ satisfies $(\mathrm T_j)$ and
    \[
        ad\equiv-j\pmod q.
    \]
    For $q\ge2$, let $\bar a$ be an inverse of $a$ modulo
    $q$; for $q=1$, put $\bar a=0$ and interpret every congruence modulo
    $1$ as vacuous. Let $X_j\subseteq I$ be the solution set of
    $(\mathrm T_j)$ found in Step~1. Then
    \begin{equation}\label{S_union}
        \mathcal S_\beta
        =
        \bigcup_{j\in\{0,1,-1,-2\}}
        \{d\in X_j\cap\mathbb Z:d\equiv-j\bar a\pmod q\}.
    \end{equation}

    Now we note that
    $\gamma t<1/t\le1/q$ on $X_0$, so
    \[
        X_0\subseteq\{t:0\le\gamma t<1/q\},
    \]
    while
    \[
    \begin{aligned}
        X_1&\subseteq\{t:\gamma t\ge1/q\},\\
        X_{-1}&\subseteq\{t:-1/q\le\gamma t<0\},\\
        X_{-2}&\subseteq\{t:\gamma t<-1/q\}
    \end{aligned}
    \] and hence the sets $X_j$ are disjoint.
    By Step~1, if $\gamma>0$, only $X_0$ and $X_1$ can be nonempty, and
    each is an interval. If $\gamma=0$, then $X_0=I$ and the other three
    sets are empty. If $\gamma<0$, only $X_{-1}$ and $X_{-2}$ can be
    nonempty. The former has at most two interval components and the
    latter has at most one. Thus at most three progression pieces occur.

    Hence, if $\gamma\ge0$, we take
    \[
        I_1:=X_0,\qquad r_1:=0,\qquad
        I_2:=X_1,\qquad r_2:=-\bar a,\qquad I_3:=\emptyset.
    \]
    If $\gamma<0$, we take
    \[
        I_1:=X_{-2},\qquad r_1:=2\bar a,\qquad r_2:=\bar a,
    \]
    and for $I_2,I_3$ the (at most two) intervals whose union is $X_{-1}$.
    In both cases \eqref{S_union} becomes \eqref{target_decomp}.
    The residues $r_1,r_2$ need not be distinct and disjointness of the
    progression pieces follows from disjointness of their supporting
    intervals.
\end{proof}

\begin{proof}[Proof of Lemma~\ref{umaze}]
    Choose $a,q$ as in \eqref{dirichlet}, and let
    $r_1,r_2$ and $I_1,I_2,I_3$ be as in Lemma~\ref{T_decomp}.
    Put $r_3:=r_2$. By \eqref{target_decomp},
    \begin{equation}\label{uma_decomp}
        \sum_{d\in I}\mu(d)\mathsf{1}_{\{\beta d\}<1/d}
        =
        \sum_{i=1}^{3}
        \sum_{\substack{d\in I_i\\d\equiv r_i\pmod q}}\mu(d).
    \end{equation}
    Expressing each congruence condition as a sum of additive characters
    gives
    \begin{equation}\label{uma_q}
        \sum_{\substack{d\in I_i\\d\equiv r_i\pmod q}}\mu(d)
        =
        \frac1q\sum_{h=0}^{q-1}
        e\!\left(-\frac{hr_i}{q}\right)
        \sum_{d\in I_i}\mu(d)e\!\left(\frac{hd}{q}\right).
    \end{equation}
    Since $I_i\subseteq[D,2D]$, Corollary~\ref{cor_zhan} and the
    triangle inequality imply
    \[
    \begin{aligned}
        \left|
        \sum_{\substack{d\in I_i\\d\equiv r_i\pmod q}}\mu(d)
        \right|
        &\le
        \frac1q\sum_{h=0}^{q-1}
        \left|\sum_{d\in I_i}\mu(d)e\!\left(\frac{hd}{q}\right)\right|\\
        &\ll_{A,\sigma}
        \frac{|I_i|}{(\log D)^A}+D^\sigma.
    \end{aligned}
    \]
    This also holds when $I_i=\emptyset$, since the corresponding sum
    vanishes. We also know that the implied constant is independent of $q$, $r_i$ and
    $I_i$ from Corollary~\ref{cor_zhan}. 

    The intervals $I_1,I_2,I_3$ are pairwise disjoint subsets of $I$, so
    \[
        \sum_{i=1}^{3}|I_i|\le |I|.
    \]
    Taking absolute values in \eqref{uma_decomp} and summing the three
    bounds therefore gives
    \[
    \begin{aligned}
        \left|
        \sum_{d\in I}\mu(d)\mathsf{1}_{\{\beta d\}<1/d}
        \right|
        &\ll_{A,\sigma}
        \sum_{i=1}^{3}
        \left(\frac{|I_i|}{(\log D)^A}+D^\sigma\right)\\
        &\ll_{A,\sigma}
        \frac{|I|}{(\log D)^A}+D^\sigma.
    \end{aligned}
    \]
    All implied constants depend only on $A$ and $\sigma$, proving
    \eqref{target_ineq} uniformly in $\beta$.
\end{proof}
\section{Large Divisor Estimates}\label{sec:large}

We now turn to the proof of Theorem~\ref{thm:main} itself. For a positive
integer $d$ and $\alpha\in\mathbb R$, we put
\[
    \mathcal A_{d,\alpha}(x)
    :=\{n\in\mathbb N:n\le x,\
    d\mid\gcd(n,\lfloor\alpha n^2\rfloor)\}.
\]
For each $n\le x$, the sum
$\sum_{d\mid\gcd(n,\lfloor\alpha n^2\rfloor)}\mu(d)$ equals $1$ if
$\gcd(n,\lfloor\alpha n^2\rfloor)=1$ and $0$ otherwise; the convention
$\gcd(n,0)=n$ covers the case $\lfloor\alpha n^2\rfloor=0$, and every
such $d$ satisfies $d\le n\le x$. Summing over $n\le x$ and
interchanging the two sums (the principle of inclusion--exclusion)
gives
\begin{equation}\label{uma_sieve}
    \mathcal S_{\alpha,2}(x)
    =\sum_{d\le x}\mu(d)|\mathcal A_{d,\alpha}(x)|.
\end{equation}
For a real parameter $R\ge1$, we split the sum in
\eqref{uma_sieve} at $d=R$:
\begin{equation}\label{uma_sieve_part}
    \begin{split}
        \mathcal S_{\alpha,2}(x)
        &=
        \sum_{d\le R}\mu(d)|\mathcal A_{d,\alpha}(x)|
        +\sum_{R<d\le x}\mu(d)|\mathcal A_{d,\alpha}(x)|\\
        &=:A_{\alpha,R}(x)+B_{\alpha,R}(x).
    \end{split}
\end{equation}
Since $\mathcal A_{d,\alpha}(x)=\varnothing$ for $d>x$, this
decomposition is valid also when $R>x$.

In this section we bound the large-divisor piece
$B_{\alpha,R}(x);$ the main goal of this section
is to prove the following theorem.

\begin{theorem}\label{thm_large}
    Let $A>0$ and $2/3<\sigma<1$. Then for all $R>3$, $x\ge1$
    and $\alpha\in\mathbb R$,
    \begin{equation}\label{uma_large}
        |B_{\alpha,R}(x)|
        \ll_{A,\sigma}
        \frac{x}{(\log R)^A}+xR^{\sigma-1},
    \end{equation}
    where the implied constant depends only on $A$ and $\sigma$.
\end{theorem}

Let us first prove a useful lemma.

\begin{lemma}\label{lem_short}
    Fix $A>0$ and $2/3<\sigma<1$. Let $D\ge2$ and let
    $I\subseteq[D,2D]$ be an interval. Then for all $x\ge1$
    and $\alpha\in\mathbb R$,
    \begin{equation}\label{uma_sieve_short}
        \left|
        \sum_{d\in I}\mu(d)|\mathcal A_{d,\alpha}(x)|
        \right|
        \ll_{A,\sigma}
        \frac{x}{(\log D)^A}+xD^{\sigma-1},
    \end{equation}
    where the implied constant depends only on $A$ and $\sigma$.
\end{lemma}

\begin{proof}
    Writing $n=dm$, we have the equivalence
    \begin{equation}\label{div_equiv}
        d\mid\lfloor\alpha m^2d^2\rfloor
        \quad\Longleftrightarrow\quad
        \{\alpha m^2d\}<1/d.
    \end{equation}
    To see this, we observe that if $t=\alpha m^2d$, then
    \[
        \lfloor dt\rfloor
        =d\lfloor t\rfloor+\lfloor d\{t\}\rfloor
    \]
    with $0\le\lfloor d\{t\}\rfloor\le d-1$; hence
    $d\mid\lfloor dt\rfloor$ if and only if
    $\lfloor d\{t\}\rfloor=0$, that is, if and only if
    $\{t\}<1/d$.

    Interchanging the finite sums therefore gives the exact identity
    \begin{equation}\label{uma_id_fin}
        \sum_{d\in I}\mu(d)|\mathcal A_{d,\alpha}(x)|
        =
        \sum_{m\le x/D}
        \sum_{\substack{d\in I\\d\le x/m}}
        \mu(d)\mathsf{1}_{\{\alpha m^2d\}<1/d}.
    \end{equation}
    For each $m$ put
    \[
        J_m:=I\cap(-\infty,x/m].
    \]
    This is an interval contained in $[D,2D]$, with $|J_m|\le D$.
    It may be empty or consist of a single point; these cases are
    included in Lemma~\ref{umaze}. Applying that lemma with
    $\beta=\alpha m^2$, and using $|J_m|\le D$, gives
    \[
        \left|
        \sum_{d\in J_m}\mu(d)
        \mathsf{1}_{\{\alpha m^2d\}<1/d}
        \right|
        \ll_{A,\sigma}
        \frac{D}{(\log D)^A}+D^\sigma,
    \]
    with an implied constant independent of $m$, $\alpha$, $x$,
    $D$, and $I$. Consequently, by \eqref{uma_id_fin} and the
    triangle inequality,
    \[
    \begin{aligned}
        \left|
        \sum_{d\in I}\mu(d)|\mathcal A_{d,\alpha}(x)|
        \right|
        &\le
        \sum_{1\le m\le x/D}
        \left|
        \sum_{d\in J_m}\mu(d)
        \mathsf{1}_{\{\alpha m^2d\}<1/d}
        \right|\\
        &\ll_{A,\sigma}
        \left\lfloor\frac{x}{D}\right\rfloor
        \left(\frac{D}{(\log D)^A}+D^\sigma\right)\\
        &\le
        \frac{x}{(\log D)^A}+xD^{\sigma-1}.
    \end{aligned}
    \]
    This also includes the case $x<D$, when the outer sum is empty,
    and proves \eqref{uma_sieve_short}.
\end{proof}

\begin{proof}[Proof of Theorem~\ref{thm_large}]
    If $x\le R$, the large-divisor sum is empty. Otherwise let
    \[
        \mathcal J
        :=\{j\in\mathbb Z_{\ge0}:2^jR<x\},
    \]
    and put
    \[
        I_j:=(2^jR,\min(2^{j+1}R,x)]
        \qquad(j\in\mathcal J).
    \]
    These pairwise disjoint intervals partition $(R,x]$.
    Applying Lemma~\ref{lem_short} with $D=2^jR$ and exponent
    $A+2$ gives
    \begin{equation}\label{ua_uma}
    \begin{split}
        |B_{\alpha,R}(x)|
        &\le
        \sum_{j\in\mathcal J}
        \left|
        \sum_{d\in I_j}\mu(d)|\mathcal A_{d,\alpha}(x)|
        \right|\\
        &\ll_{A,\sigma}
        x\sum_{j\ge0}
        \left(
        (\log R+j\log2)^{-A-2}
        +(2^jR)^{\sigma-1}
        \right)\\
        &\ll_{A,\sigma}
        \frac{x}{(\log R)^A}+xR^{\sigma-1}.
    \end{split}
    \end{equation}
    For the last step, the second series is geometric because
    $\sigma<1$, and, since $R>3$, yielding
    \[
        \sum_{j\ge0}(\log R+j\log2)^{-A-2}
        \le
        (\log R)^{-A-2}
        +(\log R)^{-A}
        \sum_{j\ge1}(j\log2)^{-2}
        \ll(\log R)^{-A}.
    \]
\end{proof}

\section{Small Divisor Estimates}\label{sec:small}

Having bounded the large-divisor sum $B_{\alpha,R}(x)$ in
Section~\ref{sec:large} using Zhan's theorem, we now turn to the
small-divisor sum $A_{\alpha,R}(x)$. For irrational $\alpha$, we
obtain quantitative discrepancy estimates for the quadratic
sequences $\{\alpha dm^2\}$ using the Erd\H{o}s--Tur\'an inequality
(Theorem~\ref{thm_et}) and the Weyl-type bound of
Lemma~\ref{uma_wooley}, both recorded in Section~\ref{sec:inputs}.

Throughout this section, $\alpha$ is irrational, $C(\alpha)$ is as in
Lemma~\ref{uma_cont_frac}, and $x\ge C(\alpha)$. We let $a=a(x)$ and
$q=q(x)$ be given by Lemma~\ref{uma_cont_frac} with $\beta=\alpha$
and $t=x$, so that
\begin{equation}\label{uma_irr_cont_*}
    \left|\alpha-\frac{a}{q}\right|\le\frac{1}{qx},
    \qquad (a,q)=1,\qquad 16\le q\le x.
\end{equation}

\begin{lemma}\label{uma_smol_div_lem}
    Let $\alpha$, $x$, $a$ and $q$ be as above. Then, uniformly for
    integers $1\le d\le q^{1/8}$,
    \begin{equation}\label{uma_smol_div}
        \left|
        |\mathcal A_{d,\alpha}(x)|-\frac{x}{d^2}
        \right|
        \ll
        \frac{x}{d^2q^{1/8}}+x^{1/2}(\log x)^{3/2},
    \end{equation}
    where the implied constant is absolute.
\end{lemma}

\begin{proof}
    For $d=1$ we have
    $|\mathcal A_{1,\alpha}(x)|=\lfloor x\rfloor$, so the left-hand
    side of \eqref{uma_smol_div} is at most $1$, and the bound holds
    because $x\ge q\ge16$. We may therefore assume that $d\ge2$.
    By \eqref{div_equiv}, writing $n=dm$,
    \[
        |\mathcal A_{d,\alpha}(x)|
        =
        \#\{m\in\mathbb N:m\le x/d,\ \{\alpha dm^2\}<1/d\}.
    \]
    Put $N_d:=\lfloor x/d\rfloor$ and
    $H:=\lfloor q^{1/4}\rfloor\ge2$.
    Since $d\le q^{1/8}\le q\le x$, we have $N_d\ge1$.
    Applying Theorem~\ref{thm_et} with $\xi_m:=\alpha dm^2$,
    $N=N_d$, $c_1=0$ and $c_2=1/d$, and using
    \[
        \frac{N_d}{d}=\frac{x}{d^2}+O\Big(\frac1d\Big),
        \qquad H+1>q^{1/4},
    \]
    we obtain
    \begin{equation}\label{uma_exp}
    \begin{split}
        \left|
        |\mathcal A_{d,\alpha}(x)|-\frac{x}{d^2}
        \right|
        &\ll\frac1d+\frac{x}{dq^{1/4}}\\
        &\quad+
        \sum_{h=1}^{H}\frac1h
        \left|\sum_{m=1}^{N_d}e(hd\alpha m^2)\right|.
    \end{split}
    \end{equation}

    Fix $1\le h\le H$. By Dirichlet's approximation theorem there
    are coprime integers $a_h,q_h$ (depending also on $d$) such that
    \begin{equation}\label{uma_diff_1}
        \left|hd\alpha-\frac{a_h}{q_h}\right|
        \le\frac{1}{2xq_h},
        \qquad 1\le q_h\le2x.
    \end{equation}
    We claim that
    \begin{equation}\label{uma_conv_ineq_1}
        q_h\ge\frac{q}{2hd}.
    \end{equation}
    Dividing \eqref{uma_diff_1} by $hd$ and comparing with
    \eqref{uma_irr_cont_*}, we get
    \[
        \left|\frac aq-\frac{a_h}{hdq_h}\right|
        \le\frac{1}{2hdxq_h}+\frac{1}{qx},
    \]
    and hence
    \[
        |ahdq_h-qa_h|
        \le\frac{q}{2x}+\frac{hdq_h}{x}
        \le\frac12+\frac{hdq_h}{x},
    \]
    since $q\le x$. If $ahdq_h\neq qa_h$, the left-hand side is a
    positive integer, so
    $1\le\frac12+\frac{hdq_h}{x}$, whence
    $q_h\ge x/(2hd)\ge q/(2hd)$. If $ahdq_h=qa_h$, then
    $q\mid hdq_h$ because $(a,q)=1$, and since $hdq_h>0$ this gives
    $hdq_h\ge q$. This proves \eqref{uma_conv_ineq_1}.
    As $h\le q^{1/4}$ and $d\le q^{1/8}$, it follows that
    \[
        q_h\ge\frac{q^{5/8}}2
        \ge\frac{16^{5/8}}2=2^{3/2}>2.
    \]
    Moreover, $q_h\le2x$ gives
    \[
        \left|hd\alpha-\frac{a_h}{q_h}\right|
        \le\frac{1}{2xq_h}\le\frac{1}{q_h^2}.
    \]
    Hence Lemma~\ref{uma_wooley} applies with $\beta=hd\alpha$,
    $r=a_h$, $s=q_h$ and $N=N_d$, and gives
    \begin{equation}\label{ineq_h_partial}
    \begin{split}
        \left|\sum_{m=1}^{N_d}e(hd\alpha m^2)\right|
        &\ll
        \frac{x}{d\sqrt{q_h}}
        +\sqrt{(x/d+q_h)\log q_h}\\
        &\ll
        \frac{x\sqrt h}{\sqrt{qd}}+\sqrt{x\log x},
    \end{split}
    \end{equation}
    where in the last step we used \eqref{uma_conv_ineq_1},
    $q_h\le2x$ and $x\ge16$.

    Inserting \eqref{ineq_h_partial} into \eqref{uma_exp}, and using
    $\sum_{h\le H}h^{-1/2}\le2H^{1/2}\le2q^{1/8}$ and
    $\sum_{h\le H}h^{-1}\ll\log q\le\log x$, we obtain
    \[
        \left|
        |\mathcal A_{d,\alpha}(x)|-\frac{x}{d^2}
        \right|
        \ll
        \frac1d+\frac{x}{dq^{1/4}}
        +\frac{x}{\sqrt d\,q^{3/8}}
        +x^{1/2}(\log x)^{3/2}.
    \]
    Finally, $1/d\le x^{1/2}(\log x)^{3/2}$, and since
    $d\le q^{1/8}$ (so that
    $d^{3/2}\le q^{3/16}\le q^{1/4}$),
    \[
        \frac{x}{dq^{1/4}}\le\frac{x}{d^2q^{1/8}},
        \qquad
        \frac{x}{\sqrt d\,q^{3/8}}\le\frac{x}{d^2q^{1/8}}.
    \]
    This proves \eqref{uma_smol_div}.
\end{proof}

\begin{lemma}\label{uma_B_lem}
    Let $\alpha$, $x$, $a$ and $q$ be as above. Then
    \begin{equation}\label{uma_B}
        A_{\alpha,q^{1/8}}(x)
        =
        \frac{x}{\zeta(2)}
        +O\left(
        \frac{x}{q^{1/8}}+x^{5/8}(\log x)^{3/2}
        \right),
    \end{equation}
    where the implied constant is absolute.
\end{lemma}

\begin{proof}
    We have
    \begin{equation}\label{B_almost_1}
    \begin{split}
        A_{\alpha,q^{1/8}}(x)
        &=\sum_{1\le d\le q^{1/8}}
        \mu(d)|\mathcal A_{d,\alpha}(x)|\\
        &=\sum_{1\le d\le q^{1/8}}\mu(d)\frac{x}{d^2}\\
        &\quad+
        \sum_{1\le d\le q^{1/8}}\mu(d)
        \left(|\mathcal A_{d,\alpha}(x)|-\frac{x}{d^2}\right)\\
        &=:\Theta_1+\Theta_2.
    \end{split}
    \end{equation}
    Since
    $\sum_{d\ge1}\mu(d)d^{-2}=\prod_p(1-p^{-2})=1/\zeta(2)$
    and $\sum_{n>y}n^{-2}\le2/y$ for $y\ge1$,
    \begin{equation}\label{theta_ineq_1}
    \begin{split}
        \Theta_1
        &=\frac{x}{\zeta(2)}
        -x\sum_{d>q^{1/8}}\frac{\mu(d)}{d^2}\\
        &=\frac{x}{\zeta(2)}
        +O\left(\frac{x}{q^{1/8}}\right).
    \end{split}
    \end{equation}
    By the triangle inequality, \eqref{uma_smol_div} and $q\le x$,
    \begin{equation}\label{theta_ineq_2}
    \begin{split}
        |\Theta_2|
        &\ll
        \sum_{1\le d\le q^{1/8}}
        \left(
        \frac{x}{d^2q^{1/8}}+x^{1/2}(\log x)^{3/2}
        \right)\\
        &\le
        \frac{\zeta(2)\,x}{q^{1/8}}
        +q^{1/8}x^{1/2}(\log x)^{3/2}\\
        &\ll
        \frac{x}{q^{1/8}}+x^{5/8}(\log x)^{3/2}.
    \end{split}
    \end{equation}
    Inserting \eqref{theta_ineq_1} and \eqref{theta_ineq_2} into
    \eqref{B_almost_1} completes the proof.
\end{proof}
\section{Proof of the Main Theorem}\label{sec:proof}
Let $\alpha\in\mathbb R\setminus\mathbb Q$. Since $q(x)$ is nondecreasing and $q(x)\to\infty$ by Lemma~\ref{uma_cont_frac}, there is $x_0(\alpha)\ge C(\alpha)$ such that $q(x)>3^8$ for all $x\ge x_0(\alpha)$. Let $x\ge x_0(\alpha)$, let $a=a(x)$ and $q=q(x)$ be as in Lemma~\ref{uma_cont_frac}, and set $R=q^{1/8}>3$. This choice yields explicit bounds for both parts of the sieve, including when $q(x)$ grows slowly.

By Lemma~\ref{uma_B_lem},
\[
    A_{\alpha,R}(x)=\frac{x}{\zeta(2)}
    +O\!\left(\frac{x}{q^{1/8}}+x^{5/8}(\log x)^{3/2}\right).
\]
For each fixed $A>0$ and $2/3<\sigma<1$, Theorem~\ref{thm_large} gives
\[
    |B_{\alpha,R}(x)|\ll_{A,\sigma}
    \frac{x}{(\log q)^A}+\frac{x}{q^{(1-\sigma)/8}},
\]
where the factor arising from $\log R=(\log q)/8$ is absorbed into
the implied constant. Thus \eqref{uma_sieve_part} yields
\begin{equation}\label{uma_final}
\begin{split}
    \mathcal S_{\alpha,2}(x)
    &=\frac{x}{\zeta(2)}\\
    &\quad+O_{A,\sigma}\!\left(
    \frac{x}{q^{1/8}}+x^{5/8}(\log x)^{3/2}
    +\frac{x}{(\log q)^A}+\frac{x}{q^{(1-\sigma)/8}}
    \right).
\end{split}
\end{equation}
For completeness, if $c=(1-\sigma)/8$, then $0<c<1/24$ and
\[
    xq^{-1/8}\le xq^{-c},\qquad
    x^{5/8}(\log x)^{3/2}\ll_\sigma x^{1-c}\le xq^{-c},
\]
using $q\le x$. Consequently the more concise estimate
\begin{equation}\label{uma_final_two_terms}
    \mathcal S_{\alpha,2}(x)=\frac{x}{\zeta(2)}
    +O_{A,\sigma}\!\left(
    \frac{x}{(\log q(x))^A}
    +\frac{x}{q(x)^{(1-\sigma)/8}}
    \right)
\end{equation}
also holds. Fixing $\sigma=3/4$ and using
$q^{-1/32}\ll_A(\log q)^{-A}$ gives \eqref{main_quantitative}.
Finally, $q(x)\to\infty$ implies that its error term is $o(x)$, proving
\eqref{conj_2} and completing the proof of Theorem~\ref{thm:main}.

\appendix
\section{Proof of Lemma~\ref{uma_wooley}}\label{app:weyl}

Whenever $\|t\|=0$, we interpret $\min(L,c/\|t\|)$ as $L$,
for $L\ge0$ and $c>0$.
We first establish the following elementary counting bound.

\begin{lemma}\label{lem_min_sum}
    Let $\beta\in\mathbb R$, and let $r,s$ be coprime integers
    with $s\ge2$ and $|\beta-r/s|\le1/s^2$.
    Then, for all integers $M,N\ge1$,
    \begin{equation}\label{min_sum}
        \sum_{m=1}^{M}
        \min\Big(N,\frac{1}{\|m\beta\|}\Big)
        \ll
        \Big(\frac Ms+1\Big)(N+s\log s),
    \end{equation}
    where the implied constant is absolute.
\end{lemma}

\begin{proof}
    We split $\{1,\dots,M\}$ into at most $M/s+1$ blocks of at
    most $s$ consecutive integers, and consider one block
    $\{m_0+j:1\le j\le s'\}$ with $1\le s'\le s$.
    Write $\beta=r/s+\delta$, so that $|\delta|\le1/s^2$.
    For $1\le j\le s'$,
    \[
        (m_0+j)\beta
        =m_0\beta+\frac{jr}{s}+j\delta,
        \qquad |j\delta|\le\frac1s,
    \]
    and hence
    \[
        \|(m_0+j)\beta\|
        \ge
        \left\|m_0\beta+\frac{jr}{s}\right\|-\frac1s.
    \]
       Since $(r,s)=1$, the points $m_0\beta+jr/s$, for $1\le j\le s'$, are
    distinct modulo $1$ and lie in $m_0\beta+\frac1s\mathbb Z$. For each
    integer $\ell\ge0$, the condition $\ell/s\le\|x\|<(\ell+1)/s$ restricts
    $x$ modulo $1$ to two half-open intervals of length $1/s$, each of
    which contains exactly one point of $m_0\beta+\frac1s\mathbb Z$. Hence
    at most two values of $j$ satisfy
    $\ell/s\le\|m_0\beta+jr/s\|<(\ell+1)/s$, and for these
    $\|(m_0+j)\beta\|\ge(\ell-1)/s$. There will be atmost $4$ terms with $\ell\le1$ and we trivially bound these by $N.$ Since $\|x\|\le1/2$, the remaining terms have
    $2\le\ell\le s/2$, so the block contributes at most
    \[
        4N+\sum_{2\le\ell\le s/2}\frac{2s}{\ell-1}\ll N+s\log s.
    \]
    Summing over the blocks proves \eqref{min_sum}.
\end{proof}

\begin{proof}[Proof of Lemma~\ref{uma_wooley}]
    For real $N\ge1$, put $L:=\lfloor N\rfloor$.
    The exponential sum is unchanged when $N$ is replaced by $L$,
    and
    \[
        \frac{L}{\sqrt{s}}+\sqrt{(L+s)\log s}
        \le
        \frac{N}{\sqrt{s}}+\sqrt{(N+s)\log s}.
    \]
    It therefore suffices to prove the estimate for integers $N\ge1$.

    Let $S:=\sum_{n=1}^{N}e(n^2\beta)$.
    Writing $n_1=n_2+h$, we have
    \[
    \begin{aligned}
        |S|^2
        &=\sum_{n_1,n_2=1}^{N}
        e\big((n_1^2-n_2^2)\beta\big)\\
        &=\sum_{|h|<N}e(h^2\beta)
        \sum_{n\in\mathcal I_h}e(2hn\beta),
    \end{aligned}
    \]
    where
    \[
        \mathcal I_h
        :=\{n\in\mathbb Z:
            1\le n\le N,\ 1\le n+h\le N\}
    \]
    is a set of at most $N$ consecutive integers.
    For every finite set $\mathcal I$ of consecutive integers
    and every $\theta\in\mathbb R$, we have
    \[
        \left|\sum_{n\in\mathcal I}e(n\theta)\right|
        \le
        \min\left(\#\mathcal I,\frac{1}{2\|\theta\|}\right).
    \]
    Since $\|2h\beta\|=\|{-2h\beta}\|$, the terms with $h$
    and $-h$ admit the same bound, and therefore
    \[
    \begin{aligned}
        |S|^2
        &\le N+2\sum_{h=1}^{N-1}
        \min\Big(N,\frac{1}{\|2h\beta\|}\Big)\\
        &\le N+2\sum_{m=1}^{2N}
        \min\Big(N,\frac{1}{\|m\beta\|}\Big).
    \end{aligned}
    \]
    By Lemma~\ref{lem_min_sum} with $M=2N$, and since
    $N\le N\log s/\log2$,
    \[
        |S|^2
        \ll
        \Big(\frac Ns+1\Big)(N+s\log s)
        \ll
        \frac{N^2}{s}+(N+s)\log s.
    \]
    Taking square roots gives \eqref{wooley}.
\end{proof}

\section*{Acknowledgements}
The author thanks his supervisor, Andrew Granville, for his invaluable guidance and suggestions throughout this work, especially for a simple presentation of Lemma~\ref{T_decomp}. This work was supported by a Master's Research Scholarship from the Fonds de recherche du Qu\'ebec (award no.~2004096;
\url{https://doi.org/10.69777/2004096}).

\section*{Declaration of generative AI use}
During the preparation of this manuscript, the author used ChatGPT
(OpenAI) and Claude (Anthropic) for assistance with exposition, LaTeX formatting, literature
searches, and the review and refinement of mathematical arguments.
The author takes full responsibility for the final manuscript,
including its mathematical content and references.

\end{document}